\documentclass[a4paper]{article}
\usepackage{amsmath,amsthm,amssymb,amscd,array}
\usepackage{fancybox}
\newtheorem{thm}{Theorem}
\newtheorem{cor}{Corollary}
\newtheorem{lem}{Lemma}

\newtheorem{property}{Property}

\begin{document}
\begin{center}
{\large A new counting methods, including the issue of 
counting labelled self-complementary graphs}
\end{center}
\begin{center}
Shinsei Tazawa\\
Department of Mathematics,
Faculty of Science and Engineering, \\ Kinki University, Higashi-Osaka
577-8502, Japan
\end{center}
abstract
\par
Harary and Palmer announced an enumeration problem of
labelled self-complementary graphs
at the end of their book (Graphical Enumeration,
Academic Press, New York and London, 1973).
This paper resolves this problem.
A method for solving this problem leads to the derivation of 
following formulas:
(a) A formula on the number of labelled graphs 
with the given order of automorphism groups of those graphs. 
(b) A formula on the number of unlabelled graphs 
with the given order of automorphism groups of those graphs.  
(c) A formula on the number of labelled self-complementary graphs 
with the given order of automorphism groups of those graphs.
(d) A formula on the number of unlabelled self-complementary graphs 
with the given order of automorphism groups of those graphs.
\\[0.5cm]
\section{Introduction}
\hspace*{5.3mm}This paper considers finite graphs.  They do not have multiple edges nor
loops on vertices.  Harary and Palmer[5] announced many graphical
enumeration problems at the end of their book, including 
an enumeration problem of labelled self-complementary
graphs.  This paper resolves the latter problem.  
A method for solving this problem leads to the derivation of 
following formulas:
(a) A formula on the number of labelled graphs 
with the given order of automorphism groups of those graphs. 
(b) A formula on the number of unlabelled graphs 
with the given order of automorphism groups of those graphs.  
(c) A formula on the number of labelled self-complementary graphs 
with the given order of automorphism groups of those graphs.
(d) A formula on the number of unlabelled self-complementary graphs 
with the given order of automorphism groups of those graphs.
\par
A formula on the number of unlabelled graphs was given by Harary[3]
by using P\'{o}lya theorem[6].  This number is found by taking the sum of
the numbers in (b)
over the orders of the automorphism groups of those graphs.
The enumeration of unlabelled self-complementary graphs dates back 
to the Read's investigation[7] in 1963 and 
Read[7] gave a formula on the number of unlabelled self-complementary graphs,
by using De Bruijn's theorem[1,2].
This number is found by taking the sum of
the numbers in (d) over the orders of 
the automorphism groups of those self-complementary graphs.
\par
In the section 2 we consider to assign a number to each labelled graph.
This assignment is very useful in the sense that 
from the number having been assigned to a graph,
the graph can be uniquely determined.
The section 3 discusses about the auotmorphism groups of graphs. 
The section 4 treats the enumerations of labelled graphs and 
unlabelled graphs with 
the orders of automorphism groups of those graphs. 
In the section 6 we treat the enumerations of labelled self-complementary
graphs, which is a main purpose in this paper, 
and the unlabelled self-complementary graphs.
\section{The indexes of graphs}
\hspace*{5.3mm}
For a positive integer $n$ let V be the set $\{1,2,...,n\}$ and let 
$\binom{V}{2}=\{\{i,j\}\subset V| i\ne j\}$. We here denote each element
$\{i,j\}$ of $\binom{V}{2}$ by $ij$ when $i<j$ and by $ji$ when $i>j$. 
For a subset
$E$ of $\binom{V}{2}$ we call the pair $(V,E)$  a labelled graph or 
a graph on $V$.  $|V|$ is called the order of the graph and $E$ is called the
edge-set of the graph.  We define a bijection $w$ from $\binom{V}{2}$
to $\{2^0=1, 2^1, 2^2,\cdots, 2^{\lambda-1}\}$ by
\begin{equation} \label{eq:wij}
 w(ij)=2^{\lambda-p_{ij}} \quad\quad ij\in\binom{V}{2},
\end{equation}
where $\lambda=\binom{n}{2}$ and $p_{ij}=\frac{(2n-i)(i-1)}{2}+(j-i)$
for $ij\in\binom{V}{2}$.
Let $\mathcal{G}_n$ be the set of all labelled graphs of 
order $n$ and let $\mathfrak{N}_n=\{0,1,2,3,\cdots,2^{\lambda}-1\}$.
We define a map $N$ from $\mathcal{G}_n$ to $\mathfrak{N}_n$
by 
\begin{equation} \label{eq:NG}
    N(G)=\left\{
     \begin{array}{ll}
       \displaystyle\sum_{ij\in E}w(ij)=
       \displaystyle\sum_{ij\in E}2^{\lambda-p_{ij}} & \quad
            \mbox{if $E$ is not empty} \\[0.2cm]
       0 & \quad \mbox{if $E$ is empty}. 
     \end{array}\right.
\end{equation}
$N(G)$ is called  the index of $G$.  We consider 
the binary number representation of an integer $L$ satisfying
$0\leqq L\leqq 2^{\lambda}-1$ and the binary number of $L$ is denoted by
$L^{(2)}$, where $L^{(2)}$ is a number of $\binom{n}{2}$ figures.  
For example a decimal number $29$ is represented as the binary number $011101$,
whose leading digit "0" is not omitted.  In the binary number $L^{(2)}$
of $L$, the $p_{ij}$-th digit is written by $L_{ij}^{(2)}$ for each
element $ij\in\binom{V}{2}$.  The binary number of the index $N(G)$ of 
a labelled graph $G$ is denoted by $N^{(2)}(G)$ and its $p_{ij}$-th
digit is written by $N_{ij}^{(2)}(G)$.
The index $N(G)$ of a labelled graph $G=(V,E)$ has 
the following property:
\begin{equation}  \label{eq:NijoneE}
N_{ij}^{(2)}(G)=1 \quad \mbox{if and only if} \quad ij\in E.
\end{equation}
We first give the following result.
\begin{thm}   \label{thm:bijection}
 The map $N$ is a bijection.
\end{thm}
\begin{proof}
Consider $L\in\mathfrak{N}_n$.  If $L=0$, then for the empty graph $G$ whose
edge-set is empty, we have $N(G)=L$.  Suppose $L\ne 0$.  Then 
the graph $G$ on $V$ having the set $E=\{ij\in\binom{V}{2} | L_{ij}^{(2)}=1\}$
as its edge-set is an element of $\mathcal{G}_n$ and it is obvious that
$N(G)=L$ holds.  For $G_1, G_2\in\mathcal{G}_n$,
if $G_1$ and $G_2$ are different, then $N^{(2)}(G_1)$ 
is not equal to $N^{(2)}(G_2)$.  This fact implies that
$N(G_1)$ is not equal to $N(G_2)$.
Hence $N$ is bijection.
\end{proof}
\hspace*{5.3mm}
As an illustration, consider a graph $G=(V,E=\{13,14,23,34\})$ on 
$V=\{1,2,3,4\}$.  Since $p_{13}=2, p_{14}=3, p_{23}=4$ and $p_{34}=6$, it is 
easy to check that $N(G)=29$ and $N^{(2)}(G)=011101$.
Conversely, a binary number $L^{(2)}=011101$ yields the labelled graph with the 
index $L$, since $L_{13}^{(2)}=L_{14}^{(2)}=L_{23}^{(2)}=L_{34}^{(2)}=1$
and $L_{12}^{(2)}=L_{24}^{(2)}=0$.
This graph is what was given in the beginning of this illustration.
\par
The complement $\overline{G}$ of a graph $G=(V,E)$ is the graph
on $V$ with edge-set $\binom{V}{2}-E$.
The index $N(\overline{G})$ of $\overline{G}$ is $2^{\lambda}-1-N(G)$. 
It is easy to see that $N^{(2)}(\overline{G})$ is obtained
by interchanging $0$ and $1$ in $N^{(2)}(G)$.
\section{Graph and group}
\hspace*{5.3mm}
Let $A$ be the symmetric group on $V$.  In this paper we assume $n\geqq 3$.
Each permutation $\alpha$ in $A$ induces a permutation $\alpha'$ 
which acts on $\binom{V}{2}$ such that for every element $ij\in\binom{V}{2}$,
\begin{displaymath}
\alpha'\{i,j\}=\{\alpha i,\alpha j\}.
\end{displaymath}
The permutation group $A'=\{\alpha' | \alpha\in A\}$ is called the pair
group of $A$.  The degree of $A'$ is $\lambda$ and $A\cong A'$, since 
$n\geqq 3$.  For a graph $G=(V,E)$ and $\alpha\in A$,
the graph which has $\{\{\alpha i,\alpha j\}\in\binom{V}{2}| ij\in E\}$ 
as its edge-set is denoted by $\alpha G$.
For a graph $G$, an element $\alpha$ in $A$ satisfying 
$\alpha G=G$ is called the automorphism of $G$.
The automorphism group $\Gamma(G)=\{\alpha\in A |\alpha G=G\}$ of $G$ is
called the group of $G$, and the number of elements which belong to the
group is called the group-order.
\par
A graph $G$ is said to be self-complementary if $G$ and $\overline{G}$ are
isomorphic, that is, if there exists $\alpha\in A$ such that 
$\alpha G=\overline{G}$ holds.  
It is well-known that every self-complementary graph $G$ has
$n\equiv 0,1 (\mod 4)$ vertices. Note that while a self-complementary
graph $G$ is isomorphic to $\overline{G}$, $G$ and $\overline{G}$
do not have the same index.
Let $\mathcal{S}_n$ be the set of all
labelled self-complementary graphs on $V$.  Then we define a set
$A_S=\{\alpha\in A | \exists G\in\mathcal{S}_n, \alpha G=\overline{G}\}$
and we also define a set $A_S(G)=\{\alpha\in A_S |\alpha G=\overline{G}\}$
for $G\in\mathcal{S}_n$.  For $\alpha'\in A'$ let $Z(\alpha')$ be the set
of all cycles in the disjoint cycle decomposition of $\alpha'$.  Then
the following lemma can easily be seen 
from the definitions of $\Gamma(G)$ and $A_S(G)$.
\begin{lem}　  \label{lem:alphaGammaGfirst}
  \begin{description}
 \item[(i)] For $\alpha\in A$ and $G=(V,E)\in\mathcal{G}_n$,
 $\alpha$ is an element of $\Gamma(G)$ if and only if 
 for each $z\in Z(\alpha')$,
    \begin{equation}  \label{eq:alphaGammaGfirst}
      ij\in E \quad\mbox{if and only if}\quad
       z(ij)\in E,\quad  \forall ij\in z
    \end{equation}
 holds.
 \item[(ii)] For $\alpha\in A$ and $G=(V,E)\in\mathcal{G}_n$,
 $G$ is a self-complementary graph satisfying $\alpha G=\overline{G}$,
 that is, $\alpha\in A_S(G)$ if and only if 
 for each $z\in Z(\alpha')$,
     \begin{equation} \label{eq:alphaASGfirst}
        ij\in E \quad\mbox{if and only if}\quad
         z(ij)\not\in E, \quad \forall ij\in z
    \end{equation}
 holds.
 \end{description}
\end{lem}
\begin{lem}　\label{lem:alphaGammaGsecond}
\begin{description}
 \item[(i)] For $\alpha\in A$ and $G\in\mathcal{G}_n$,
 $\alpha$ is an element of $\Gamma(G)$ if and only if
 for each $z\in Z(\alpha')$,
     \begin{equation} \label{eq:alphaGammaGsecond}
       N_{ij}^{(2)}(G)=1\,\, \mbox{{\rm for}}\,\, \forall ij\in z\quad
       \mbox{or} \quad
       N_{ij}^{(2)}(G)=0 \,\, \mbox{{\rm for}}\,\,\forall  ij\in z
     \end{equation}
  holds.
 \item[(ii)] For $\alpha\in A$ and $G\in\mathcal{G}_n$,
 $G$ is a self-complementary graph satisfying $\alpha G=\overline{G}$,
 that is, $\alpha\in A_S(G)$  if and only if
 for each $z\in Z(\alpha')$,
     \begin{equation} \label{eq:alphaASGsecond}
        N_{ij}^{(2)}(G)+N_{z(ij)}^{(2)}(G)=1
         \quad \mbox{{\rm for}}\,\, \forall ij\in z
     \end{equation}
 holds.
 \end{description}
\end{lem}
\begin{proof}
(i) is immediate from the combination of (\ref{eq:NijoneE}) and
(i) in Lemma \ref{lem:alphaGammaGfirst}.  We next prove (ii).
Let $E$ be the edge-set of $G$.  Then
(\ref{eq:NijoneE}) can be written as follows:
\begin{equation} \label{eq:zijnotENij}
 z(ij)\not\in E \quad\mbox{if and only if}\quad
 N_{z(ij)}^{(2)}(G)=0.
\end{equation}
Therefore, it is easy to see from (\ref{eq:NijoneE}),
(\ref{eq:zijnotENij}) and (ii) in Lemma \ref{lem:alphaGammaGfirst}
 that the desired result is obtained.
\end{proof}
\begin{thm} \label{thm:property}
Let $\mathcal{K}(P)$ be the set of all labelled graphs of order $n$ 
all of which have a property $P$ and let $\xi$ be a positive integer.
If the group-order of the group of every graph in $\mathcal{K}(P)$ is $\xi$,
then the number of graphs 
which belong to $\mathcal{K}(P)$ and which
are not mutually isomorphic is ,that is, the number of unlabelled graphs 
in $\mathcal{K}(P)$ is given by
\begin{equation} \label{eq:kpn}
 \frac{\xi|\mathcal{K}(P)|}{n!}.
\end{equation}
\end{thm}
\begin{proof}
Let $l(G)$ be the number of ways of labelling a graph $G$.
Then using a well-known $l(G)=\frac{n!}{|\Gamma(G)|}\,\,$[5;p4],
we get the equality (\ref{eq:kpn}), since $|\Gamma(G)|=\xi$ for 
$G\in\mathcal{K}(P)$. This complete the proof.
\end{proof}
\hspace*{5.3mm}
We have the following results with respect to self-complementary
graphs.
\begin{thm} \label{thm:ASGGammaG}
Suppose that $G$ is an element of $\mathcal{S}_n$. Then
for any $\alpha\in A_S(G)$, $\alpha^{-1} A_S(G)=\Gamma(G)$ holds.
\end{thm}
\begin{proof}
It follows that for any $\beta\in A_S(G)$, 
$(\alpha^{-1}\beta)G=\alpha^{-1}(\beta G)=\alpha^{-1}\overline{G}=G$ holds.
Therefore we have $\alpha^{-1}\beta\in\Gamma(G)$, which implies 
$\alpha^{-1}A_S(G)\subseteq\Gamma(G)$.  Conversely,
Any element $\gamma$ of $\Gamma(G)$ can be written as 
$\gamma=\alpha^{-1}\alpha\gamma$.  Thus 
$(\alpha\gamma)G=\alpha(\gamma G)=\alpha G=\overline{G}$ is obtained.
Therefore, we have $\alpha\gamma\in A_S(G)$, that is, 
$\gamma\in\alpha^{-1}A_S(G)$, which implies
$\Gamma(G)\subseteq\alpha^{-1}A_S(G)$.  Hence 
$\alpha^{-1} A_S(G)=\Gamma(G)$ holds.
\end{proof}
\hspace*{5.3mm}
Theorem\ref{thm:ASGGammaG} gives the following.
\begin{cor} \label{cor:ASGGammaGcardinality}
For $G\in\mathcal{S}_n$, $|A_S(G)|=|\Gamma(G)|$ holds.
\end{cor}
For $\alpha\in A$ let $m_k(\alpha)$ be the number of cycles of length $k$
in the disjoint cycle decomposition of $\alpha$.  Then we prove
\begin{thm} \label{lem:m-characterize}
Let $n$ be a positive integer satisfying $n\equiv 0,1(\bmod\enskip 4)$.
Then $\alpha$ is an element of $A_S$ if and only if
the following condition holds:
\begin{eqnarray}
    m_1(\alpha)&&\le1\quad\mbox{and} \label{eq:jone} \\[0.1cm]
m_k(\alpha)&&\left\{
      \begin{array}{rl}
          &\ge1\quad\mbox{if $k\equiv 0\pmod{4}$}  \label{eq:fourteenth}\\
          &=0\quad\mbox{if $k\not\equiv 0\pmod{4}$}
          \end{array}\right.
\end{eqnarray}
for $k=2,3,...,n$.
\end{thm}
\begin{proof}
Suppose $\alpha\in A_S$.
Then there exists $G\in\mathcal{S}_n$ satisfying $\alpha G=\overline{G}$.
Accordingly, it is easy to see that the length of each cycle 
in $Z(\alpha')$ must
be even.  This fact implies $m_1(\alpha)\le1$ and $m_k(\alpha)=0$ for 
odd $k(\geqq 3)$.  If $\alpha$ has a cycle of length $2k$ providing
k is odd, then it generates a cycle of odd length $k$ in $\alpha'$
which contradicts the fact that the length of each cycle in $Z(\alpha')$
is even.
\par
On the other hand, suppose that the permutation $\alpha$ of $A$ satisfies
the conditions (\ref{eq:jone}) and (\ref{eq:fourteenth}).
Then the length of each cycle in $Z(\alpha')$ for the permutation
$\alpha'\in A'$ induced by $\alpha$ is always even.  Accordingly,
we can construct the binary number $L^{(2)}$ of $\binom{n}{2}$ figures 
satisfying 
\begin{displaymath}
 L_{ij}^{(2)}+L_{z(ij)}^{(2)}=1, \quad \forall ij\in z
\end{displaymath}
for $z\in Z(\alpha')$.
If we consider the labelled graph $G$ on $V$ which has 
$E=\{ij\in\binom{V}{2} | L_{ij}^{(2)}=1\}$ as its edge-set, 
(ii) in Lemma \ref{lem:alphaGammaGsecond} tells us that $G$ is 
a self-complementary graph of order $n$ satisfying $\alpha G=\overline{G}$.
Hence we have $\alpha\in A_S$.
\end{proof}
\section{Counting of graphs}
\hspace*{5.3mm}
For each $z\in Z(\alpha')$ in $\alpha\in A$, we consider
\begin{equation} \label{eq:Wzsumw}
 W(z)=\sum_{ij\in z}w(ij).
\end{equation}
Then to $\alpha\in A$ there corresponds the following polynomial of $x$:
\begin{equation} \label{eq:Falphanx}
 F_{\alpha}^{(n)}(x)=\prod_{z\in Z(\alpha')}f_z^{(n)}(x),
\end{equation}
where
\begin{equation} \label{eq:fznx}
 f_z^{(n)}(x)=1+x^{W(z)}.
\end{equation}
Let us consider the generating function for labelled graphs of order $n$
\begin{equation} \label{eq:Fnx}
  F^{(n)}(x)=\sum_{\alpha\in A}F_{\alpha}^{(n)}(x).
\end{equation}
Let $c_L^{(\alpha)}$ be the coefficient of $x^L$ in the power series
expansion of $F_{\alpha}^{(n)}(x)$
for $\alpha\in A$,
and let $c_L$ be 
the coefficient of $x^L$ in the power series
expansion of $F^{(n)}(x)$.  Then the following
equality is an immediate consequence of (\ref{eq:Fnx}):
\begin{equation} \label{eq:cncnalpha}
c_L=\sum_{\alpha\in A}c_L^{(\alpha)}.
\end{equation}
Let $f(x)$ be a polynomial of $x$.  Then by the symbol $a_nx^n\vDash f(x)$,
$a_nx^n$ will mean the term in the polynomial $f(x)$.  On the other hand,
If $a_nx^n$ is not any term of $f(x)$, that is, if $a_nx^n\nvDash f(x)$,
we use the convention $a_n=0$.
\par
Let $lg(n)=|\mathcal{G}_n|$. Then $lg(n)$ is equal to $2^{\lambda}$, 
since $|\mathcal{G}_n|=2^{\lambda}$, as well-known.
Consider an $n!$ by $lg(n)$ matrix $M(\mathcal{G}_n)=||c_L^{(\alpha)}||$,
which has a row for each $\alpha\in A$ and a column for 
each $G\in\mathcal{G}_n$, provided that $c_L^{(\alpha)}=0$ when
$c_L^{(\alpha)}x^L\nvDash F_{\alpha}^{(n)}(x)$.
Since by Theorem \ref{thm:bijection} there is a one-to-one correspondence
between $\mathcal{G}_n$ and $\mathfrak{N}_n$, the column corresponding
to $G\in\mathcal{G}_n$ is referred to as the $N(G)$-column.
\begin{property} \label{property:functionpropertyfirst}
In the matrix $M(\mathcal{G}_n)$, if $c_L^{(\alpha)}\ne 0$, then
$c_L^{(\alpha)}=1$ holds.
\end{property}
\begin{proof}
Consider $\alpha\in A$.  
If $L=0$, then since the constant term in
the power series expansion of $F_{\alpha}^{(n)}(x)$ is $1$, 
we have $c_0^{(\alpha)}=1$.
Suppose that $L\ne 0$.  Then note first that 
$c_L^{(\alpha)}x^L\vDash F_{\alpha}^{(n)}(x)$, since $c_L^{(\alpha)}\ne 0$.
From the observation of (\ref{eq:Falphanx})
there exists a subset $K=\{z_1,z_2,\cdots,z_t\}$ of $Z(\alpha')$
such that $L=W(z_1)+W(z_2)+\cdots+W(z_t)$.
It can easily seen
that the binary number $L^{(2)}$ satisfies
\begin{eqnarray*}
  & & L_{ij}^{(2)}=1,\,\,\forall ij\in z, \forall z\in K \\
  & & L_{ij}^{(2)}=0,\,\,\forall ij\in z, \forall z\in Z(\alpha')-K.
\end{eqnarray*}
Thus if we consider the labelled graph $G\in\mathcal{G}_n$ whose edge-set is 
$\{ij\in\binom{V}{2} | L_{ij}^{(2)}=1\}$, then it follows that
its index $N(G)$ is $L$ and that $\alpha G=G$ holds 
by (i) in Lemma \ref{lem:alphaGammaGsecond}.
Since $L^{(2)}$ is determined by $L$, vice versa, 
$L$ is determined by $L^{(2)}$,
the subset $K$ being discussed just above can be uniquely determined by $L$.
This implies that $c_L^{(\alpha)}=1$.
\end{proof}
\hspace*{5.3mm}
This property states that for $\alpha\in A$, every coefficient
of the power series expansion of $F_{\alpha}^{(n)}(x)$ is $1$.
\begin{property} \label{property:functionpropertysecond}
The coefficient $c_L$ of $x^L$ in $F^{(n)}(x)$ is equal to 
the group-order of the group of labelled graph having the index $L$.
\end{property}
\begin{proof}
In the matrix $M(\mathcal{G}_n)$, (\ref{eq:cncnalpha}) denotes the $L$-th
column sum.  The $L$-th column is associated with
the labelled graph $G$ having the index $N(G)=L$, and 
(\ref{eq:cncnalpha}) becomes 
$c_L=\displaystyle\sum_{\alpha\in\Gamma(G)}c_L^{(\alpha)}$.  
Applying Property \ref{property:functionpropertyfirst}, we get
$c_L=|\Gamma(G)|$.
\end{proof}
\hspace*{5.3mm}
If we consider two sets: 
$\mathcal{O}(\mathcal{G}_n)=\{|\Gamma(G)|| G\in\mathcal{G}_n\}$
and $\mathcal{C}(\mathcal{G}_n)=\{c_L | c_Lx^L\vDash F^{(n)}\}$, then
Property \ref{property:functionpropertysecond} gives the
following theorem, which states that the list of
the group-orders of the groups of
all graphs of order $n$ can be obtained by expanding $F^{(n)}(x)$.
\begin{thm} \label{thm:labelledorderd}
$\mathcal{O}(\mathcal{G}_n)=\mathcal{C}(\mathcal{G}_n)$
holds.
\end{thm}
For $\xi\in\mathcal{O}(\mathcal{G}_n)$, let $lg(n;\xi)$ be
the number of {\it labelled} graphs $G$ of order $n$ such that
$|\Gamma(G)|=\xi$, and let $g_{\xi}^{(n)}$ be the number
of {\it unlabelled} graphs $G$ of order $n$ such that
$|\Gamma(G)|=\xi$.
For a positive integer $\xi$, put
$\mathcal{L}_{\xi}(\mathcal{G}_n)=
\{L | c_L=\xi,\,c_Lx^L\vDash F^{(n)}(x)\}$.
Then we have the following three theorems.
\begin{thm} \label{thm:lgnxiLxi}
The equality
\begin{equation} \label{eq:lgnxiLxi}
 lg(n;\xi)=|\mathcal{L}_{\xi}(\mathcal{G}_n)|, \quad \xi\in\mathcal{C}(\mathcal{G}_n)
\end{equation}
holds.
\end{thm}
\begin{proof}
As seen from the statement of Property \ref{property:functionpropertysecond} and
its proof, 
each labelled graph $G$ in $\mathcal{G}_n$ having $|\Gamma(G)|=\xi$ 
is associated with the $N(G)$-th column of $M(\mathcal{G}_n)$
whose column sum is $\xi$, and 
each element of $\mathcal{L}_{\xi}(\mathcal{G}_n)$
is also associated with a column of $M(\mathcal{G}_n)$
whose column sum is $\xi$.
Accordingly, $lg(n;\xi)$ is the number of columns
of $M(\mathcal{G}_n)$ whose column sum is $\xi$
and also $|\mathcal{L}_{\xi}(\mathcal{G}_n)|$ is 
the number of columns
of $M(\mathcal{G}_n)$ whose column sum is $\xi$.
Hence we obtain the equality (\ref{eq:lgnxiLxi}),
noting Theorem \ref{thm:labelledorderd}.
\end{proof}
\begin{thm}　\label{thm:unlabelledgxi}
 The equality
 \begin{equation} \label{eq:unlabelledorderd}
   g_{\xi}^{(n)}=\frac{\xi}{n!}|\mathcal{L}_{\xi}(\mathcal{G}_n)|,\,\,\,\,\quad
 \xi\in\mathcal{C}(\mathcal{G}_n)
 \end{equation}
 holds.
\end{thm}
\begin{proof}
In Theorem \ref{thm:property}, let $P$ denote the property
that the group-order of the group of a labelled graph in 
    $\mathcal{G}_n$ is $\xi$.
    Then $\mathcal{K}(P)$ is the set of all labelled graphs $G$ in 
    $\mathcal{G}_n$ satisfying $|\Gamma(G)|=\xi$, that is, 
    $|\mathcal{K}(P)|=lg(n;\xi)$.
    Hence using Theorem \ref{thm:property}, 
    Theorem \ref{thm:labelledorderd} and Theorem \ref{thm:lgnxiLxi},
    the equality (\ref{eq:unlabelledorderd}) is obtained.
\end{proof}
\hspace*{5.3mm}
Let $g(n)$ be the number of unlabelled graphs of order $n$.  A formula 
on $g(n)$ has been given by Harary[3].  But an enumeration method
having been here presented will derive a formula on $g(n)$, which 
is stated in the following.
\begin{thm}
 The number $g(n)$ of unlabelled graphs of order $n$ is given by
 \begin{equation} \label{eq:unlabelledgn}
    g(n)=
    \sum_{\xi\in \mathcal{C}(\mathcal{G}_n)}\frac{\xi}{n!}
    |\mathcal{L}_{\xi}(\mathcal{G}_n)|.
 \end{equation}
\end{thm}
\begin{proof}
  Since $g(n)=
    \displaystyle\sum_{\xi\in\mathcal{O}(\mathcal{G}_n)}g_{\xi}^{(n)}$,
    using Theorem \ref{thm:labelledorderd} and 
    Theorem \ref{thm:unlabelledgxi},
    we get the equality (\ref{eq:unlabelledgn}).
\end{proof}
\section{Numerical results for graphs}
\hspace*{5.3mm}
We have the numerical examples for $n\geqq 3$, since $n\geqq 3$ in this 
paper.  Table 1 gives the numerical examples for labelled graphs.
Each pair $(\xi, lg(n;\xi))$ in the Table 1 shows the number $lg(n;\xi)$ 
of labelled
graphs of order $n$ with the group-order $\xi$ 
for a few small values of $n\geqq 3$. The last column of Table 1 
shows the number $lg(n)$ of labelled graphs of order $n$, which is well-known.
On the other hand,
Table 2 gives the numerical examples for unlabelled graphs.
Each pair $(\xi, g_{\xi}^{(n)})$ in the Table 2 shows the number
$g_{\xi}^{(n)}$
of unlabelled
graphs of order $n$ with the group-order $\xi$.  The last column of Table 2 
shows the number $g(n)$ of unlabelled graphs of order $n$.
\begin{table}[h]
\begin{center}
\caption{Labelled graphs}
\begin{displaymath}
\renewcommand{\arraystretch}{2.5}
 \begin{array}{c|l|l} \hline
 n   & \multicolumn{1}{c|}{(\xi,\,lg(n;\xi)\,)} & lg(n) \\ \hline
 3       &  (2,6),(6,2) & 2^3 \\ \hline
 4   &  (2,36),(4,12),(6,8),(8,6),(24,2) & 2^6  \\ \hline
 5   &  
  \begin{minipage}[c]{10cm}
    (2,660),(4,180),(6,40),(8,60),(10,12), (12,60), (24,10),(120,2) 
  \end{minipage} & 2^{10}  \\[2mm] \hline
 6  
  & 
  \begin{minipage}[c]{10cm}
(1,5760), (2,16560), (4,6480), (6,960), (8,1260), (10,144), (12,1080),
(16,270), (24,60), (36,40), (48,120), (72,20), (120,12), (720,2)
\end{minipage} & 2^{15} \\[2mm] \hline
7  
 &
 \begin{minipage}[c]{10cm}
 (1,766080), (2,892080), (4,312480), (6,31920), (8,46620), 
 (10,1008), (12,29400), (14,720), (16,6300), (20,1008), (24,5040), 
 (36,840), (48,2940), (72,280), (120,84), (144,210), (240,126), 
 (720,14), (5040,2)
 \end{minipage} & 2^{21} \\[2mm] \hline
 \end{array}
\end{displaymath}
\end{center}
\end{table}
\vspace*{0.5cm}
\begin{table}[h]
\begin{center}
\caption{Unlabelled graphs}
\begin{displaymath}
\renewcommand{\arraystretch}{2.5}
  \begin{array}{c|l|l} \hline
  n  & \multicolumn{1}{c|}{(\xi,\,g_{\xi}^{(n)}\,)} & g(n) \\ \hline
  3  & (2,2), (6,2) & 4 \\ \hline
  4  &  (2,3), (4,2), (6,2), (8,2), (24,2) & 11 \\ \hline 
 5  & (2,11), (4,6), (6,2), (8,4), (10,1), (12,6), (24,2), (120,2) & 34 
    \\[2mm] \hline
  6  & \begin{minipage}[c]{10cm}
      (1,8), (2,46), (4,36), (6,8), (8,14), (10,2), (12,18), (16,6),
       (24,2), (36,2), (48,8), (72,2), (120,2), (720,2)
     \end{minipage} 
     & 156  \\[2mm] \hline 
  7  & \begin{minipage}[c]{10cm}
      (1,152), (2,354), (4,248), (6,38), (8,74), (10,2), (12,70), 
     (14,2), (16,20), (20,4), (24,24), (36,6), (48,28), (72,4), 
     (120,2), (144,6), (240,6), (720,2), (5040,2)
     \end{minipage} 
    & 1044  \\[2mm] \hline 
     \end{array}
\end{displaymath}
\end{center}
\end{table}
\newpage
\section{Counting of self-complementary graphs}
\hspace*{5.3mm}
In this section we make a consideration about the enumeration of
self-complementary graphs.  This enumeration is made along the argument
in the preceding section.
Let $n$, throughout this section, be an integer satisfying 
$n\equiv 0,1 (\mod 4)$ and $ n\geqq 4$.  For $\alpha\in A_S$, 
the length of each cycle in $Z(\alpha')$ is always even.  
We denote a cycle $z\in Z(\alpha')$ by 
$z=(i_1j_1,i_2j_2,\cdots,i_{k-1}j_{k-1},i_kj_k)$.
The cycle can be represented in some forms in the usual cyclic
representation, as well-known.  For example, the cycle$(12,23,34,14)$ can 
also be written 
$(23,34,14,12),\,\,(34,14,12,23)$ or $(14,12,23,34)$.  We assume that
the leading position of $z$ is occupied by the element at which the value of
$w$ is the maximum among those elements of $z$.  Let 
$P_1(z)=\{i_1j_1,i_3j_3,\cdots,i_{k-1}j_{k-1}\}$ be the set of elements
in the odd position in the cycle $z$ and 
let $P_2(z)=\{i_2j_2,i_4j_4,\cdots,i_kj_k\}$
be the set of elements in the even position in $z$.  Then, 
note that for $ij\in z$, $ij\in P_1(z)$ if and only if $z(ij)\in P_2(z)$.
For each $z\in Z(\alpha')$ in $\alpha\in A_S$, we consider
\begin{equation} \label{eq:woddwevensum}
  W_1(z)=\sum_{ij\in P_1(z)}w(ij) \quad \mbox{and} \quad 
  W_2(z)=\sum_{ij\in P_2(z)}w(ij).
\end{equation}
Then to $\alpha\in A_S$ there corresponds the following polynomial of $x$:
\begin{equation} \label{eq:Salphanx}
 U_{\alpha}^{(n)}(x)=\prod_{z\in Z(\alpha')}u_z^{(n)}(x),
\end{equation}
where
\begin{equation}  \label{eq:eq:sznx}
 u_z^{(n)}(x)=x^{W_1(z)}+x^{W_2(z)}, \quad z\in Z(\alpha').
\end{equation}
Let us consider 
the generating function for labelled self-complementary graphs of
order $n$
\begin{equation} \label{eq:Snx}
  U^{(n)}(x)=\sum_{\alpha\in A_S}U_{\alpha}^{(n)}(x).
\end{equation}
Let $lsc(n)=|\mathcal{S}_n|$ be the number of labelled self-complementary
graphs of order $n$. Then we shall show that $lsc(n)$ is equal to
the number of terms in the power series expansion of $U^{(n)}(x)$.
For $G\in\mathcal{S}_n$, there exists $\alpha\in A_S$
satisfying $\alpha G=\overline{G}$.  It is seen from (ii) in 
Lemma \ref{lem:alphaGammaGsecond} that the index $N(G)$ of $G$
satisfies 
\begin{equation} \label{eq:GLij}
 N_{ij}^{(2)}(G)+N_{z(ij)}^{(2)}(G)=1, \quad \forall ij\in z
\end{equation}
for each $z\in Z(\alpha')$.  For $z\in Z(\alpha')$, if 
we put $P(z)=\{ij\in z|N_{ij}^{(2)}(G)=1\}$ , then it follows from 
(\ref{eq:GLij}) that $P(z)$ agrees with either $P_1(z)$ or $P_2(z)$. 
Furthermore, if we put $W(z)=\displaystyle\sum_{ij\in P(z)}w(ij)$
for each $z\in Z(\alpha')$, then 
it is easy to check that $N(G)=\displaystyle\sum_{z\in Z(\alpha')}W(z)$.  Thus
$\displaystyle\prod_{z\in Z(\alpha')}x^{W(z)}=x^{N(G)}$
is a term in the power series expansion of $U_{\alpha}^{(n)}(x)$,
which implies that this term is appeared in the power series expansion
of $U^{(n)}$.  Since $N$ is a bijection from $\mathcal{G}_n$ to
$\mathfrak{N}_n$ by Theorem \ref{thm:bijection},
we have the following theorem. 
\begin{thm} \label{thm:labelledselfconplementary}
Let $\mathcal{L}(\mathcal{S}_n)=\{L |d_Lx^L\vDash U^{(n)}(x)\}$.
Then the equality $lsc(n)=|\mathcal{L}(\mathcal{S}_n)|$ holds.
\end{thm}
Let $d_L^{(\alpha)}$ be the coefficient of $x^L$ in the power series
expansion of $U_{\alpha}^{(n)}(x)$
for $\alpha\in A_S$, and let $d_L$ be the coefficient of $x^L$ in 
the power series expansion of $U^{(n)}(x)$.
Then the following equality is an immediate consequence of (\ref{eq:Snx}):
\begin{equation} \label{eq:selfcncnalpha}
d_L=\sum_{\alpha\in A_S}d_L^{(\alpha)}.
\end{equation}
 Consider an $|A_S|$ by $lsc(n)$ matrix $M(\mathcal{S}_n)
=||d_L^{(\alpha)}||$, which has a row for each $\alpha\in A_S$ and
a column for each $G\in\mathcal{S}_n$, provided that $d_L^{(\alpha)}=0$
when $d_L^{(\alpha)}x^L\nvDash U_{\alpha}^{(n)}(x)$.  $M(\mathcal{S}_n)$
may be regarded as a submatrix of $M(\mathcal{G}_n)$, if $c_L^{(\alpha)}$
in the latter is replaced by $d_L^{(\alpha)}$ in the former.
Thus the column corresponding to $G\in\mathcal{S}_n$ is referred to
as the $N(G)$-column.
\begin{property} \label{property:selffunctionpropertyfirst}
In the matrix $M(\mathcal{S}_n)$, if $d_L^{(\alpha)}\ne 0$, then
$d_L^{(\alpha)}=1$ holds.
\end{property}
\begin{proof}
Consider $\alpha\in A_S$. Note first that
$d_L^{(\alpha)}x^L\vDash U_{\alpha}^{(n)}(x)$,
since $d_L^{(\alpha)}\ne 0$.  Also, note that $L\ne 0$, 
since we assume that $n\geqq 4$.  Put $Z(\alpha')=\{z_1,z_2,\cdots,z_k\}$
and put $I=\{1,2\}$.  Then
from the observation of (\ref{eq:Salphanx}) 
there exists a $k-$tuple 
$(l_1,l_2,\cdots,l_k)\in\underbrace{I\times I\times\cdots\times I}_{k}$
such that $L=W_{l_1}(z_1)+W_{l_2}(z_2)+\cdots+W_{l_k}(z_k)$.  
It can easily seen that
the binary number $L^{(2)}$ satisfies
\begin{eqnarray*}
  & & L_{ij}^{(2)}=1,\,\,\forall ij\in
  P_{l_1}(z_1)\cup P_{l_2}(z_2)\cup\cdots\cup P_{l_k}(z_k) \\
  & & L_{ij}^{(2)}=0,\,\,\forall ij\in 
  P_{l'_1}(z_1)\cup P_{l'_2}(z_2)\cup\cdots\cup P_{l'_k}(z_k),
\end{eqnarray*}
where $l'_h\ne l_h,\,\, l'_h\in I\,(h=1,2,\cdots,k)$.
Thus if we consider the labelled 
graph $G\in\mathcal{G}_n$ whose edge-set is 
$\{ij\in\binom{V}{2} | L_{ij}^{(2)}=1\}$,
then it is obvious that its index $N(G)$ is $L$, and
(ii) in Lemma \ref{lem:alphaGammaGsecond} tells us that $G$ is
a self-complementary graph satisfying $\alpha\in A_S(G)$.
As seen at the end of the proof of
Property \ref{property:functionpropertyfirst},
since $L^{(2)}$ is determined by $L$, vice versa, 
$L$ is determined by $L^{(2)}$,
the $k$-tuple $(l_1,l_2,\cdots,l_k)$
being discussed just above can be uniquely determined by $L$.
This implies that $d_L^{(\alpha)}=1$.
\end{proof}
\hspace*{5.3mm}
This property states that for $\alpha\in A_S$, every coefficient
of the power series expansion of $U_{\alpha}^{(n)}(x)$ is $1$.
\begin{property} \label{property:selffunctionpropertysecond}
The coefficient $d_L$ of $x^L$ in $U^{(n)}(x)$ is equal to 
the group-order of the group of labelled self-complementary 
graph having the index $L$.
\end{property}
\begin{proof}
In the matrix $M(\mathcal{S}_n)$, (\ref{eq:selfcncnalpha}) denotes the $L$-th
column sum.  The $L$-th column is associated with
the labelled self-complementary graph $G$ having the index $N(G)=L$, and 
(\ref{eq:selfcncnalpha}) becomes 
$d_L=\displaystyle\sum_{\alpha\in A_S(G)}d_L^{(\alpha)}$.  
Applying Property \ref{property:selffunctionpropertyfirst}, we get
$d_L=|A_S(G)|$.  Hence we have $d_L=|\Gamma(G)|$
from Corollary \ref{cor:ASGGammaGcardinality}.
\end{proof}
\hspace*{5.3mm}
If we consider two sets: 
$\mathcal{O}(\mathcal{S}_n)=\{|\Gamma(G)|| G\in\mathcal{S}_n\}$
and $\mathcal{D}(\mathcal{S}_n)=\{d_L | d_Lx^L\vDash U^{(n)}\}$, then
Property \ref{property:selffunctionpropertysecond} gives the
following theorem, which states that the list of the group-orders of the groups
of all self-complementary graphs of order $n$ can be obtained by
expanding $U^{(n)}(x)$.
\begin{thm} \label{thm:labelledselforderd}
$\mathcal{O}(\mathcal{S}_n)=\mathcal{D}(\mathcal{S}_n)$
holds.
\end{thm}
For $\xi\in\mathcal{O}(\mathcal{S}_n)$, let $lsc(n;\xi)$ be
the number of {\it labelled} self-complementary graphs $G$ 
of order $n$ such that
$|\Gamma(G)|=\xi$, and let $s_{\xi}^{(n)}$ be the number
of {\it unlabelled} self-complementary graphs $G$ of order $n$ such that
$|\Gamma(G)|=\xi$.
For a positive integer $\xi$, put
$\mathcal{L}_{\xi}(\mathcal{S}_n)=
\{L | d_L=\xi,\,d_Lx^L\vDash U^{(n)}(x)\}$.
Then we have the following three theorems.
\begin{thm} \label{thm:lscxiLxi}
The equality
\begin{equation} \label{eq:iscxiLxi}
 lsc(n;\xi)=|\mathcal{L}_{\xi}(\mathcal{S}_n)|, 
 \quad \xi\in\mathcal{D}(\mathcal{S}_n)
\end{equation}
holds.
\end{thm}
\begin{proof}
As seen from the statement of
Property \ref{property:selffunctionpropertysecond}
and its proof,
each labelled graph $G$ in $\mathcal{S}_n$ having $|\Gamma(G)|=\xi$ 
is associated with the $N(G)$-th column of $M(\mathcal{S}_n)$
whose column sum is $\xi$, and
each element of $\mathcal{L}_{\xi}(\mathcal{S}_n)$
is also associated with a column of $M(\mathcal{S}_n)$
whose column sum is $\xi$.
Accordingly, $lsc(n;\xi)$ is the number of columns
of $M(\mathcal{S}_n)$ whose column sum is $\xi$,
and also $|\mathcal{L}_{\xi}(\mathcal{S}_n)|$ is 
the number of columns
of $M(\mathcal{S}_n)$ whose column sum is $\xi$.
Hence we obtain the equality (\ref{eq:iscxiLxi}),
noting Theorem \ref{thm:labelledselforderd}.
\end{proof}
\begin{thm}　\label{thm:unlabelledselfxi}
 The equality
 \begin{equation} \label{eq:unlabelledselforderd}
   s_{\xi}^{(n)}=\frac{\xi}{n!}|\mathcal{L}_{\xi}(\mathcal{S}_n)|,\,\,\,\,\quad
 \xi\in\mathcal{D}(\mathcal{S}_n)
 \end{equation}
 holds.
\end{thm}
\begin{proof}
In Theorem \ref{thm:property}, let $P$ denote the property
that the group-order of the group of a labelled self-complementary graph in 
    $\mathcal{S}_n$ is $\xi$.
    Then $\mathcal{K}(P)$ is the set of all labelled 
    self-complementary graphs $G$ in 
    $\mathcal{S}_n$ satisfying $|\Gamma(G)|=\xi$, that is, 
    $|\mathcal{K}(P)|=lsc(n;\xi)$.
    Hence using Theorem \ref{thm:property}, 
    Theorem \ref{thm:labelledselforderd} and Theorem \ref{thm:lscxiLxi},
    the equality (\ref{eq:unlabelledselforderd}) is obtained.
\end{proof}
\hspace*{5.3mm}
Let $sc(n)$ be the number of unlabelled self-complementary graphs of order $n$.
A formula on $sc(n)$ has been given by Read[7].  But an enumeration method
having been here presented will derive a formula on $sc(n)$, which 
is stated in the following.
\begin{thm}
 The number $sc(n)$ of unlabelled self-complementary graphs of order $n$ is 
 given by
 \begin{equation} \label{eq:unlabelledselfscn}
    sc(n)=
    \sum_{\xi\in \mathcal{D}(\mathcal{S}_n)}\frac{\xi}{n!}
    |\mathcal{L}_{\xi}(\mathcal{S}_n)|.
 \end{equation}
\end{thm}
\begin{proof}
  Since $sc(n)=
    \displaystyle\sum_{\xi\in\mathcal{O}(\mathcal{S}_n)}s_{\xi}^{(n)}$,
    using Theorem \ref{thm:labelledselforderd} and 
    Theorem \ref{thm:unlabelledselfxi},
    we get the equality (\ref{eq:unlabelledselfscn}).
\end{proof}
\section{Numerical results for self-complementary graphs}
\hspace*{5.3mm}
This section gives the numerical results for self-complementary graphs
of a few small orders $n$ satisfying $n\equiv 0,1 (\mod 4)$. 
Table 3 gives the numerical examples for labelled self-complementary graphs.
Each pair $(\xi,\,lsc(n;\xi))$ in the Table 3 shows the number $lsc(n;\xi)$ 
of labelled self-complementary graphs of order $n$ with the group-order $\xi$. 
The last column of Table 3 
shows the number $lsc(n)$ of labelled self-complementary graphs of order $n$.
On the other hand,
Table 4 gives the numerical examples for unlabelled self-complementary graphs.
Each pair $(\xi,\,s_{\xi}^{(n)})$ in the Table 4 
shows the number $s_{\xi}^{(n)}$ of unlabelled self-complementary
graphs of order $n$ with the group-order $\xi$.  The last column of Table 4 
shows the number $sc(n)$ of unlabelled self-complementary graphs of order $n$.
\begin{table}[h]
\begin{center}
\caption{Labelled self-complementary graphs}
\begin{displaymath}
\renewcommand{\arraystretch}{1.5}
 \begin{array}{c|l|l} \hline
   n  & \multicolumn{1}{c|}{(\xi,\,lsc(n;\xi)\,)}  & lsc(n) \\ \hline
   4  &  (2,12)                               & 12            \\ \hline
   5  &  (2,60),(10,12)                        & 72     \\ \hline
   8  &  (2,60480),(4,20160),(8,15120),(32,2520) & 98280 \\ \hline
   9  &  
   \begin{minipage}[c]{9cm}
     (2,3265920),(4,544320),(8,226800), (20,36288),\\ (32,45360),
     (72,5040) 
   \end{minipage}
                       & 4123728 \\[2mm] \hline
  \end{array}
\end{displaymath}
\end{center}
\end{table}
\vspace*{0.5cm}
\begin{table}[h]
\begin{center}
\caption{Unlabelled self-complementary graphs}
\begin{displaymath}
\renewcommand{\arraystretch}{1.5}
 \begin{array}{c|l|l} \hline
    n  &  \multicolumn{1}{c|}{(\xi,\,s_{\xi}^{(n)}\,)} & sc(n) \\ \hline
 4  & (2,1) & 1 \\ \hline
 5  & (2,1), (10,1)  & 2 \\ \hline
 8  & (2,3), (4,2), (8,3), (32,2)   & 10 \\ \hline
 9  &  (2,18), (4,6), (8,5), (20,2), (32,4), (72,1)
      & 36 \\ \hline
 \end{array}
\end{displaymath}
\end{center}
\end{table}
\\
\vspace*{0.5cm}
Acknowledgements
\\
\hspace*{5.3mm}
The author wishes to express his deepest gratitude  to
Professor Fumikazu Tamari, who is in 
Fukuoka University of Education, for frequent, stimulating, and helpful
discussions, and for his kind suggestions with respect to  the format
of this paper.  The author also wishes to thank Professors Yasuo Ohno and 
Tsunenobu Asai, who are in Department of Mathematics, Kinki University,
for their helpful advice.
\\[2cm]
References\\
\begin{description}
\item[[1\!\!]] N.G. de Bruijn, Generalization of P\'{o}lya's
fundamental theorem in enumeration combinatorial analysis, Indagationes
Math. 21(1959), 59-69.
\item[[2\!\!]] N.G. de Bruijn, P\'{o}lya's theory of counting,in:
E. F. Beckenbach(Ed.), Applied Combinatorial Mathematics:
Wiley, New York, 1964, pp.144-184.
\item[[3\!\!]] F. Harary, The number of linear, directed, rooted, and
connected graphs, Trans. Amer. Math. Soc. 78(1955), 445-463.
\item[[4\!\!]] F. Harary, Graph Theory, 
Addison-Wesley Pub., Massachusetts,1972.
\item[[5\!\!]] F. Harary and E.M. Palmer, Graphical Enumeration,
Academic Press, New York and London, 1973.
\item[[6\!\!]] G. Polya, Kombinatorische Anzahlbestimmungen f\"{u}r Gruppen,
Graphen und chemische Verbindungen, Acta Math. 68(1937), 145-254.
\item[[7\!\!]] R.C. Read, On the number of self-complementary graphs
and digraphs, Journal London Math. Soc.  {\bf 38}(1963), 99-104.
\end{description}
\end{document}